\documentclass[11pt]{article}
\usepackage{latexsym}
\usepackage{amsfonts}

\newcommand{\be}{\begin{equation}}
\newcommand{\ee}{\end{equation}}
\newcommand{\bea}{\begin{eqnarray}}
\newcommand{\eea}{\end{eqnarray}}
\newcommand{\bean}{\begin{eqnarray*}}
\newcommand{\eean}{\end{eqnarray*}}
\newcommand{\brray}{\begin{array}}
\newcommand{\erray}{\end{array}}

\newcommand{\newsection}[1]{\setcounter{equation}{0} \setcounter{dfn}{0}
\section{#1}}

\newtheorem{dfn}{Definition}[section]
\newtheorem{thm}[dfn]{Theorem}
\newtheorem{lmma}[dfn]{Lemma}
\newtheorem{ppsn}[dfn]{Proposition}
\newtheorem{crlre}[dfn]{Corollary}
\newtheorem{xmpl}[dfn]{Example}
\newtheorem{rmrk}[dfn]{Remark}

\newcommand{\bdfn}{\begin{dfn}}
\newcommand{\bthm}{\begin{thm}}
\newcommand{\blmma}{\begin{lmma}}
\newcommand{\bppsn}{\begin{ppsn}}
\newcommand{\bcrlre}{\begin{crlre}}
\newcommand{\bxmpl}{\begin{xmpl}}
\newcommand{\brmrk}{\begin{rmrk}}

\newcommand{\edfn}{\end{dfn}}
\newcommand{\ethm}{\end{thm}}
\newcommand{\elmma}{\end{lmma}}
\newcommand{\eppsn}{\end{ppsn}}
\newcommand{\ecrlre}{\end{crlre}}
\newcommand{\exmpl}{\end{xmpl}}
\newcommand{\ermrk}{\end{rmrk}}

\newcommand{\IC}{\mathbb{C}}

\newcommand{\eps}{\epsilon}
\newcommand{\LMD}{\Lambda}
\newcommand{\cla}{{\cal A}}
\newcommand{\clb}{{\cal B}}
\newcommand{\clh}{{\cal H}}
\newcommand{\clk}{{\cal K}}
\newcommand{\cll}{{\cal L}}

\def \bbs {\mbox{\boldmath $s$}}
\def \bbt {\mbox{\boldmath $t$}}

\newcommand{\prf}{\noindent{\it Proof\/}: }
\newcommand{\ots}{\otimes}
\newcommand{\raro}{\rightarrow}
\newcommand{\RARO}{\Rightarrow}
\newcommand{\seq}{\subseteq}
\newcommand{\ol}{\overline}
\newcommand{\lgl}{\langle}
\newcommand{\rgl}{\rangle}
\newcommand{\one}{1\!\!1}

\newcommand{\NI}{\noindent}
\newcommand {\CC}{\centerline}
\def \qed { \mbox{}\hfill $\Box$\vspace{1ex}}

\newcommand{\hta}{\hat{\cla}}
\newcommand{\whtG}{\widehat{G}}
\newcommand{\kernel}{\mbox{ker\,}}
\newcommand{\ran}{\mbox{ran\,}}
\newcommand{\tlds}{\tilde{S}}
\begin{document}

\author{{\large Arupkumar Pal}\\
         Indian Statistical Institute,\\[-.5ex]
         7, SJSS Marg, New Delhi--110\,016, INDIA\\[-.5ex]
         email: arup@isid.ac.in}
\title{Regular Operators on Hilbert $C^*$-modules}
\date{September\,11, 1997\\ Revised: April\,3,1998}
\maketitle
\vspace{-2ex}

   \begin{quotation}
\CC{\bf Abstract}
A regular operator $T$ on a Hilbert $C^*$-module is defined
just like a closed operator on a Hilbert space, with the extra
condition that the range of $(I+T^*T)$ is dense. Semiregular
operators are a slightly larger class of operators that may not
have this property. It is shown that, like in the case of regular
operators, one can, without any loss in generality, restrict
oneself to semiregular operators on $C^*$-algebras. We then prove
that for abelian $C^*$-algebras as well as for
subalgebras of the algebra of compact operators, any closed
semiregular operator is automatically regular. We also determine
how a regular operator and its extensions (and restrictions) are
related. Finally, using these results, we give a criterion for
a semiregular operator on a liminal $C^*$-algebra to have a
regular extension.\\[1ex]
{\bf AMS Subject Classification No.:} {\large 46}H{\large 25},
                                  {\large 47}C{\large 15}
   \end{quotation}

\newsection{Introduction}\label{s1}
Hilbert $C^*$-modules were first studied by Kaplansky (\cite{kap})
for abelian $C^*$-algebras, and later for more general $C^*$-algebras by
Rieffel (\cite{rie}) and Paschke (\cite{pa}).
Kasparov developed the theory further, and used them
to get deep and far-reaching results (\cite{k}) in KK-theory.
As the name suggests,
Hilbert $C^*$-modules are very similar to Hilbert spaces,
with $C^*$-algebra elements playing the role of scalars. With
the development of quantum groups and noncommutative geometry,
the study of Hilbert $C^*$-modules
has assumed further importance. In noncommutative geometry,
for example, the role of vector bundles on (noncommutative) spaces is
played by Hilbert $C^*$-modules. For a locally compact quantum
group $G$, the $C^*$-algebra $\cla=C_0(G)$ of
`continuous vanishing-at-infinity functions'
on $G$ is a Hilbert $C^*$-module over itself.
Analogues of bounded operators in the Hilbert $C^*$-module
context are {\em adjointable operators}. For example, for any
$C^*$-algebra considered as a Hilbert $C^*$-module over itself,
adjointable operators are the elements of its multiplier algebra.
In various contexts that the Hilbert $C^*$-modules arise, one also needs to
study `unbounded adjointable operators', or what are now known
as {\em regular operators}. These were first introduced by
Baaj \& Julg in \cite{bj}, where they gave a nice construction
of Kasparov bimodules in KK-theory using regular operators.
Later they were rediscovered by Woronowicz (\cite{wo1})
while investigating noncompact quantum groups. He considered
$C^*$-algebras rather than general Hilbert $C^*$-modules
(we shall see in section~\ref{s3} that there is no loss in
generality in doing this), and called them elements affiliated to
the $C^*$-algebra. `Coordinate functions' on locally compact noncompact
quantum groups are examples of such objects. Representations of
locally compact noncompact groups (quantum as well as classical)
are examples of regular operators on more general Hilbert $C^*$-modules.
Lance gave a brief indication in his book (\cite{la}) about the
possible role Hilbert modules  might play in
studying representations of quantum groups.

Let us quickly recall the definition  of a regular operator.
Let $\cla$ be a $C^*$-algebra. An operator $T$  from a Hilbert $\cla$-module
$E$ to another Hilbert $\cla$-module $F$ is said to be
{\em regular} if \\[-3ex]
\begin{enumerate}\renewcommand{\theenumi}{\alph{enumi}}
\item $T$ is closed and densely defined,
\item its adjoint $T^*$ is also densely defined, and
\item range of $I+T^*T$ is dense in $F$.
\end{enumerate}
Note that if we set $\cla = \IC$, i.e.\ if we take $E$ and $F$ to
be Hilbert spaces, then this is exactly the definition of a
closed operator, except that in that case, both the second and
the third condition follow from the first one. In the
Hilbert $C^*$-module context, one needs to add these extra conditions
in order to get a reasonably good theory. But when one deals
with specific unbounded operators on concrete Hilbert $C^*$-modules,
it is usually extremely difficult to
verify the last condition, though the first two conditions are
relatively easy to check. So it would be interesting to find other
more easily manageable conditions that are equivalent to the
last condition above. In \cite{wo1}, Woronowicz
gave a criterion based on the graph of an operator for it to be
regular, and to this date, this remains the only attempt
in this direction.

In the present paper, we will consider a somewhat larger class of
operators that we call semiregular operators, which are,
roughly speaking, operators satisfying the first two conditions above.
We then investigate the following two problems, namely, under what
conditions are they regular, and when do they admit regular extensions.
We will assume elements of $C^*$-algebra theory
as can be found for example in Dixmier (\cite{dix}). For an account
on the $C^*$-module theory required, we refer the reader to
Lance~(\cite{la}).

\vspace{1ex}

 {\bf Notations}. $\clh$ will denote a complex separable Hilbert space.
$\clb_{0}(\clh)$ and $\clb(\clh)$ will respectively stand for
the space of compact operators on $\clh$ and the space of all bounded
operators on it. $\cla$  denotes a $C^*$-algebra, usually nonunital,
and $M(\cla)$ is its multiplier algebra. The algebra $\cla$ will
always be assumed to be separable.
$\pi$, with or without sub- (or super-) scripts
will usually denote representations of the $C^*$-algebra under consideration,
and $\clh_\pi$ will be the Hilbert space on which the representation acts.
$E$ and  $F$ will denote Hilbert $C^*$-modules, generally
over a $C^*$-algebra $\cla$. $\lgl E,E\rgl$ will denote
the linear span of $\{\lgl x,y\rgl:x,y\in E\}$ in $\cla$.
For a Hilbert $C^*$-module $E$, $\clk(E)$ and $\cll(E)$ will
denote respectively the space of all `compact' operators on $E$ and
the space of all adjointable operators on $E$.
$S$, $T$, $\bbs$, $\bbt$ etc.\ will be operators on
Hilbert modules or on Hilbert spaces. For an operator $T$,
$G(T)$ will denote its graph.
For a topological space $X$, $C_0(X)$ (respectively $C_c(X)$)
will denote the algebra of continuous functions on $X$
vanishing at infinity (resp.\ with compact support).
\vspace{1ex}

Before we end this section, let us state here a Stone-Weirstrass type
theorem for $C^*$-algebras that will be very useful in studying
regular operators on Hilbert $C^*$-modules.
\bthm \label{t1.1}
Let $\cla$ be a separable $C^*$-algebra, $\hta$ being its spectrum.
Let $J$ be a right ideal in $\cla$ such that $\pi(J)$ is dense
in $\pi(\cla)$ for all $\pi\in\hta$. Then $J$ is dense in $\cla$.
\ethm

Normally, in Stone-Weirstrass type results, the subspace $J$ is assumed
to be a $*$-subalgebra, and the proof goes roughly like this:
if $J$ is not dense in $\cla$, one constructs a nonzero state on $\cla$
that vanishes on $J$. The corresponding GNS representation must
also vanish on $J$. Using separability, one can now get a point
$\pi\in\hta$  that vanishes on $J$. Since this is
not the case, one reaches a contradiction. Thus the key step in the
proof is the construction of the state, which is made possible by
the condition that $J$ is closed under involution. In our case, $J$ is
not necessarily $*$-closed; but as lemma~2.9.4 in \cite{dix} tells us,
the condition that it is a right ideal is strong enough to guarantee
that such a construction is still possible.


\newsection{Semiregular operators} \label{s2}
Let us start with the following definition.

\bdfn \label{t2.1}
Let $E$ and $F$ be Hilbert $\cla$-modules. An operator $T:E\raro F$ is called
{\em semiregular} if \\[-3ex]
\begin{enumerate}\renewcommand{\theenumi}{\alph{enumi}}
\item $D_T$ is a dense right submodule in $E$ {\rm (i.e.\ $D_T\cla\seq D_T$)},
\item $T$ is closable,
\item $T^*$ is densely defined.
\end{enumerate}
\edfn

Observe that from (c), it follows that $T$ is $\cla$-linear.
Any regular operator is of course semiregular. But there are also
many semiregular operators that are not regular, as
the following example illustrates.

Let   $\cla = C[0,1]$, and $E=C[0,1]\ots L_2(0,1)$. Let
\[ \brray{lcl}
D &=& \{f\in L_2(0,1):\; f \mbox{ absolutely continuous, }
           f'\in L_2(0,1)\},\\
D_0 &=& \{f\in L_2(0,1):\; f \mbox{ absolutely continuous, }
           f'\in L_2(0,1), f(0)=f(1) \},\\
D_{00} &=& \{f\in L_2(0,1):\; f \mbox{ absolutely continuous, }
           f'\in L_2(0,1), f(0)=f(1)=0 \}.
\erray
\]
Define $T$ on $D$ by $Tf=if'$. Let $T_0=T|_{D_0}$. Now define an operator
$\bbt$ on $E$ as follows:
\[
D(\bbt)=\{f\in E: f_0\in D_{00}, f_\pi\in D_0 \mbox{ for } 0<\pi\leq 1,
             \pi\mapsto {f_\pi}' \mbox{ continuous}\},
\]
\[
(\bbt f)(\pi)=i{f_\pi}'.
\]

\bppsn \label{t2.2}
The operator $\bbt$ defined above is a closed semiregular nonregular
operator.
\eppsn
\prf Let us first of all show that
\[
D(\bbt)_\pi=\cases{D_{00} & if $\pi=0$,\cr
                   D_0 & if $0<\pi\leq 1$.}
\]
{}From the definition of $D(\bbt)$, it is clear that
$D(\bbt)_0\seq D_{00}$ and $D(\bbt)_\pi\seq D_0$ for $0<\pi\leq 1$.
To show the reverse inclusions, take $f\in D_{00}$.
Define $g(\pi,x)=f(x)$. Check that $g\in D(\bbt)$. Therefore
$D(\bbt)_0=D_{00}$.
Next, fix some $\pi_0\in (0,1]$ and take $f\in D_0$.
This time, take $g(\pi,x)=\frac{\pi}{\pi_0}f(x)$.
Then $g_{\pi_0}=f$ and  $g\in D(\bbt)$, so that $D(\bbt)_{\pi_0}=D_0$.

It is easy to check that for $f, g\in D(\bbt)$,
$\lgl \bbt f,g\rgl=\lgl f,\bbt g\rgl$. Therefore $\bbt\seq\bbt^*$.
{}From this and from the fact that
$C[0,1]\ots_{\mbox{\footnotesize alg}}D_{00}\seq D(\bbt)$,
it is clear that $\bbt$ is semiregular. To show that it is closed, take
$f_n\in D(\bbt)$  such that $f_n$ converges to $f$ and $\bbt f_n$
converges to $\hat{f}$. Then $(f_n)_\pi\in D(\bbt)_\pi$,
$f_\pi=\lim (f_n)_\pi$ and ${\hat{f}}_\pi=\lim (\bbt f_n)_\pi$.
Since $(f_n)_\pi\in D(\bbt)_\pi$, $(\bbt f_n)_\pi=i(f_n)_\pi'$.
Therefore by closedness of $\bbt_\pi$, $f_\pi\in D(\bbt)_\pi$,
$i{f_\pi}'={\hat{f}}_\pi$. $\pi\mapsto i{f_\pi}'={\hat{f}}_\pi$
is continuous and hence $f\in D(\bbt)$, $\hat{f}=\bbt f$.

Finally, if $\bbt$ is regular, then
$\{(I+\bbt^*\bbt)f:f\in D(\bbt^*\bbt)\}$
is dense in $E$, so that for any $\pi$, $\{(I+\bbt^*\bbt)f(\pi):f\in
D(\bbt^*\bbt)\}$
is dense in $E_\pi=L_2(0,1)$. Now
\[
\{(I+\bbt^*\bbt)f(0):f\in D(\bbt^*\bbt)\}
  \seq \{f-f'': f\in D_{00}, f'\in D_0\}.
\]
Notice that the right hand side = $\ran (I+T_0T^*)$. So
its orthogonal complement is given by
$\kernel(I+TT_0^*)=\kernel(I+TT_0)
=\{f:f\in D_0, f'\in D, f=f''\}$. Hence the function
$f:x\mapsto \exp(x)+\exp(1-x)\in\kernel(I+TT_0)$.
Therefore $\{(I+\bbt^*\bbt)f(0):f\in D(\bbt^*\bbt)\}$ can not be dense
in $L_2(0,1)$, which means, $\bbt$ can not be regular.
\qed

This operator $\bbt$ shares many features with regular operators.
Here are two of them.
\bppsn \label{t2.3}
The operator $\bbt$ satisfies the following:\\
{\rm 1.} $\bbt^*$ is regular,\\
{\rm 2.} $D(\bbt^*\bbt)$ is a core for $\bbt$.
\eppsn
\prf
1. We will first show that $D(\bbt^*)_\pi=D_0$ for all $\pi$.
Take any $f\in D(\bbt^*)$. Then for any $g\in D(\bbt)$,
$\lgl \bbt g$, $f\rgl=\lgl g$, $\bbt^* f\rgl$. So
$\lgl \bbt_\pi g_\pi$, $f_\pi\rgl=\lgl g_\pi$, $(\bbt^* f)_\pi\rgl$.
This means $f_\pi\in D({\bbt_\pi}^*)$ and ${\bbt_\pi}^*f_\pi=(\bbt^*f)_\pi$.
In our context, $f_0\in D$, $f_\pi\in D_0$ for $0<\pi\leq 1$. Therefore
$D(\bbt^*)_0\seq D$, $D(\bbt^*)_\pi\seq D_0$ for $0<\pi\leq 1$.
We have already abserved that $\bbt\seq\bbt^*$. Hence
$D(\bbt^*)_\pi=D_0$ for $0<\pi\leq 1$. Now choose an $f\in D(\bbt^*)$.
Then for any $\pi\in(0,1]$,
\[
\brray{rclcl}
\lgl \one,\bbt^*f\rgl(\pi)&=& i(f_\pi(1)-f_\pi(0))&=&0,\\
\lgl \one,\bbt^*f\rgl(0)&=& i(f_0(1)-f_0(0)).&&
\erray
\]
By continuity, $f_0(1)=f_0(0)$, i.e.\ $D(\bbt^*)_0\seq D_0$.
To show the reverse inclusion, take $f\in D_0$. Define
$\tilde{f}(\pi,x)=f(x)$, $h(\pi,x)=if'(x)$. Then for any $g\in D(\bbt)$,
$\lgl \bbt g$, $\tilde{f}\rgl=\lgl g, h\rgl(\pi)$. Therefore
$\tilde{f}\in D(\bbt^*)$. Since ${\tilde{f}}_0=f$, we get
$D_0\seq D(\bbt^*)_0$.

Let us denote $\bbt^*$ by $\bbs$ in this paragraph.
{}From what we have seen, ${\bbs_\pi}^*=\bbs_\pi$ for all $\pi$.
Hence $D({\bbs_\pi}^*\bbs_\pi)=\{f\in D_0: f'\in D_0\}$.
Clearly $D(\bbs^*\bbs)_\pi\seq D({\bbs_\pi}^*\bbs_\pi)$.
Take any $f\in D({\bbs_\pi}^*\bbs_\pi)$. Define
$\tilde{f}(\pi,x)=f(x)$, $h(\pi,x)=-f''(x)$. Then for any
$g\in D(\bbs)$, $\lgl \bbs g, \bbs\tilde{f}\rgl=\lgl g,h\rgl$.
Therefore $\tilde{f}\in D(\bbs^*\bbs)$ and $\bbs^*\bbs \tilde{f}=h$.
Since ${\tilde{f}}_\pi=f$, we have
$D(\bbs^*\bbs)_\pi=D({\bbs_\pi}^*\bbs_\pi)$. This implies
$\{(I+\bbs^*\bbs)g(\pi):g\in D(\bbs^*\bbs)\}
  = \{(I+{\bbs_\pi}^*\bbs_\pi)f:f\in D({\bbs_\pi}^*\bbs_\pi)\}.$
Each $\bbs_\pi$ is closed, so the right hand side is dense in $L_2(0,1)$.
Consequently,$\{(I+\bbs^*\bbs)g:g\in D(\bbs^*\bbs)\}$ is dense
in $E$ (essentially by theorem~\ref{t1.1}), which means $\bbs$
is regular.\vspace{1ex}

2. It is easy to see that
\[
D(\bbt^*\bbt)_\pi\seq\cases{\{f\in D_0:f'\in D_0\}& if $0<\pi\leq 1$,\cr
                             \{f\in D_{00}:f'\in D_0\}& if $\pi=0$.}
\]
For any $f\in D_{00}$ such that $f'\in D_0$, the function
$g(\pi,x)=f(x)$ can easily be seen to be in $D(\bbt^*\bbt)$. Therefore
in the second case above, we actually have equality. To show that
equality holds in the second case as well, choose an $f\in D_0$
for which $f'\in D_0$ and define $g(\pi,x)=\frac{\pi}{\pi_0}f(x)$.
One can then verify that $g\in D(\bbt^*\bbt)$. So we now have
\[
D(\bbt^*\bbt)_\pi=\cases{\{f\in D_0:f'\in D_0\}& if $0<\pi\leq 1$,\cr
                             \{f\in D_{00}:f'\in D_0\}& if $\pi=0$.}
\]

We now prove that $D(\bbt^*\bbt)_\pi$ is a core for $\bbt_\pi$
for all $\pi$. Notice that for $\pi>0$, $D(\bbt^*\bbt)_\pi=D({T_0}^2)$,
and $\bbt_\pi=T_0$. So the above assertion holds in this case.
Suppose $\pi=0$. Take an $f\in D_{00}=D(\bbt_0)=D(T|_{D_{00}})$.
Choose $g_n\in D_{00}$ such that
\[
{g_n}'(0)=-f'(0),\quad {g_n}'(1)=-f'(1),
        \quad \|g_n\|\raro 0,\quad \|{g_n}'\|\raro 0,
\]
then $f+g_n\in D_{00}$, $(f+g_n)'\in D_{00}\seq D_0$, and
\[
\|f+g_n-f\|=\|g_n\|\raro 0,\quad \|(f+g_n)'-f'\|=\|{g_n}'\|\raro 0.
\]
Thus $D(\bbt^*\bbt)_0$ is a core for $\bbt_0$.

We are now ready to show that $D(\bbt^*\bbt)$ is a core for $\bbt$.
Take an $f\in D(\bbt)$. Let $\eps>0$ be any given number.
Choose a partition
$0=\pi_0<\pi_1<\ldots <\pi_n=1$ of $[0,1]$ such that
\[
\|f_\pi-f_{\pi'}\|<\eps,\quad \|{f_\pi}'-{f_{\pi'}}'\|<\eps,
\]
whenever $\pi$ and $\pi'$ belong to the same subinterval.
Choose $h_i\in D(\bbt^*\bbt)_{\pi_i}$ satisfying
\[
\|f_{\pi_i}-h_i\|<\eps,\quad \|{f_{\pi_i}}'-{h_i}'\|<\eps.
\]
Define $g$ as follows:
\[
g(\pi,x)=\frac{\pi_{i+1}-\pi}{\pi_{i+1}-\pi_i}h_i(x)
            + \frac{\pi-\pi_i}{\pi_{i+1}-\pi_i}h_{i+1}(x),
\quad \pi_i\leq\pi\leq \pi_{i+1}.
\]
Then
$g_{\pi_i}=h_i$, ${g_\pi}'\in D_0$ for all $\pi$,
$g_\pi\in D_0$ for all $\pi>0$, $g_0\in D_{00}$ and the maps
$\pi\mapsto {g_\pi}'$ and $\pi\mapsto{g_\pi}''$ are continuous.
{}From these, one can now easily check that
$g\in D(\bbt^*\bbt)$, and
\[
\|g_\pi-f_\pi\|<5\eps, \quad \|{g_\pi}'-{f_\pi}'\|<5\eps.
\]
Hence we have $\|g-f\|=\sup_\pi\|g_\pi-f_\pi\|<5\eps$
and $\|\bbt g-\bbt f\|=\sup_\pi\|{g_\pi}'-{f_\pi}'\|<5\eps$.
Thus $D(\bbt^*\bbt)$ is a core for $\bbt$.
\qed

\NI{\bf Remarks:} 1. The example constructed above is very similar
in spirit to an example of a nonregular selfadjoint operator
first constructed by Hilsum~(\cite{hil}).

2. Propositions \ref{t2.2} and \ref{t2.3} together imply that
for a closed semiregular operator, regularity of its
adjoint does not ensure regularity of the original operator,
i.e.\ corollary~9.6 in \cite{la} is false.

\newsection{Semiregular operators on $C^*$-algebras} \label{s3}
Let $\cla$ be a $C^*$-algebra and $E$ be a Hilbert $\cla$-module.
Any regular operator on $E$ is
uniquely determined by its $z$-transform (or the bounded transform)
which is an element $z$ of $\cll(E)$ satisfying the following two conditions:
\be \label{e3.1}
\|z\|\leq 1,\quad (1-z^*z)^{1/2}E \mbox{ dense in } E.
\ee
The space $\clk(E)$ of compact operators on $E$ is a $C^*$-algebra,
and hence is  a Hilbert $C^*$-module over itself. Any regular
operator on this $C^*$-algebra is uniquely determined by its
$z$-transform $w$ in $\cll(\clk(E))$ that obeys the following:
\be \label{e3.2}
\|w\|\leq 1,\quad (1-w^*w)^{1/2}\clk(E) \mbox{ dense in } \clk(E).
\ee
Now, via the isomorphism $\cll(E)\cong\cll(\clk(E))$, the set of
elements in $\cll(E)$ satisfying (\ref{e3.1}) can be identified with the
set of elements in $\cll(\clk(E))$ satisfying (\ref{e3.2}), which means,
regular operators on $E$ can be identified with regular operators on $\clk(E)$.
One can therefore deal with regular operators on $C^*$-algebras without any
loss of generality. In this section, we will show that the same is true for
semiregular operators as well. The main trouble in this case is that we
do not have $z$-transforms at our disposal any more.

Let us denote by $R(E)$ the space of all regular operators on a
Hilbert $C^*$-module $E$, and by $SR(E)$ the space of all semiregular
operators on $E$. Define two maps $\phi_1:SR(E)\raro SR(\clk(E))$
and
$\phi_2:SR(\clk(E))\raro SR(E)$
as follows: let $S\in SR(\clk(E))$, $T\in SR(E)$. Let
\bean
D(\phi_1(T))&:=&\mbox{span}\,\{|x\rgl\lgl y|:x\in D_T, y\in E\},\\
\phi_1(T)|x\rgl\lgl y| &:=&
     |Tx\rgl\lgl y|,\quad |x\rgl\lgl y|\in D(\phi_1(T));\\
D(\phi_2(S)) &:=& \mbox{span}\,\{ax:a\in D_S, x\in E\},\\
\phi_2(S)(ax) &:=& (Sa)x,\quad ax\in D(\phi_2(S)).
\eean
{}From the semiregularity of $S$ and $T$, it follows that $\phi_1(T)$ and
$\phi_2(S)$ are well defined and semiregular. Let us now list
some properties of these two maps $\phi_1$ and $\phi_2$.

\blmma\label{t3.1}
Let $S, S_1, S_2\in SR(\clk(E))$ and $T, T_1, T_2\in SR(E)$.
Then we have
\[
\brray{rcl}
T_1\seq T_2 & \RARO & \phi_1(T_1)\seq \phi_1(T_2),\\
S_1\seq S_2 & \RARO & \phi_2(S_1)\seq \phi_2(S_2)
\erray
\]
\elmma
\prf Straightforward.\qed

\blmma\label{t3.2}
Let $S\in SR(\clk(E))$, $T\in SR(E)$. Then
\[
\phi_1(T)\seq \phi_1(\ol{T})\seq \ol{\phi_1(T)},
\]
\[
\phi_2(S)\seq \phi_2(\ol{S})\seq \ol{\phi_2(S)}.
\]
\elmma
\prf The first inclusion in both cases follows from part~2 of the previous
lemma. For the second inclusion, take
$\sum|x_i\rgl\lgl y_i|\in D(\phi_1(\ol{T}))$,
where the $x_i$'s come from $D_{\ol{T}}$. There exist $x_i^{(n)}\in D_T$
such that $x_i=\lim x_i^{(n)}$ and $\ol{T}x_i=\lim_n Tx_i^{(n)}$.
Therefore $\sum|x_i\rgl\lgl y_i|=\lim_n \sum|x_i^{(n)}\rgl\lgl y_i|$,
and $\phi_1(\ol{T})(\sum|x_i\rgl\lgl y_i|)=\sum|\ol{T}x_i\rgl\lgl y_i|
            =\lim_n \phi_1(T)(\sum|x_i^{(n)}\rgl\lgl y_i|)$.
So
$\bigl(\sum|x_i\rgl\lgl y_i|, \phi_1(\ol{T})(\sum|x_i\rgl\lgl y_i|)\bigr)
      \in G(\ol{\phi_1(T)})$.
Thus $G(\phi_1(\ol{T}))\seq G(\ol{\phi_1(T)})$,
i.e.\ $\phi_1(\ol{T})\seq \ol{\phi_1(T)}$.

Next, take $(ax, \phi_2(\ol{S})(ax))\in G(\phi_2(\ol{S}))$. Then
$a\in D_{\ol{S}}$, so that there exist $a_n\in D_S$ such that
$a=\lim_n a_n$ and $\ol{S} a=\lim_n Sa_n$. Therefore
$ax=\lim_n a_nx$, and $(\ol{S} a)x=\lim_n (Sa_n)x$. Since
$a_nx\in D_{\phi_2(S)}$ and $\phi_2(S)(a_nx)=(Sa_n)x$, it follows that
$(ax,(\ol{S} a)x)\in G(\ol{\phi_2(S)})$. Therefore
$\phi_2(\ol{S})\seq\ol{\phi_2(S)}$.\qed

\blmma \label{t3.3}
Let $S\in SR(\clk(E))$, $T\in SR(E)$. Then
\[
\brray{rll}
1.&\phi_1(\phi_2(S))\seq S, & \phi_1(\ol{\phi_2(S)})\seq \ol{S},\\
2.&\phi_2(\phi_1(T))\seq T, & \phi_2(\ol{\phi_1(T)})\seq \ol{T}.
\erray
\]
\elmma
\prf The first inclusion in both cases is trivial. The second inclusion
follows from the first and the forgoing lemma.\qed

\blmma \label{t3.4}
Let $S\in SR(\clk(E))$, $T\in SR(E)$. Let $E$ be separable. Then\\
{\rm 1.} $D(\phi_1(\phi_2(S)))$ is a core for $\ol{S}$,\\
{\rm 2.} $D(\phi_2(\phi_1(T)))$ is a core for $\ol{T}$.
\elmma
\prf By definition,
\bean
D(\phi_1(\phi_2(S)))&=&\mbox{span}\,\{|y\rgl\lgl z|: y\in D(\phi_2(S)), z\in
E\}\\
&=& \mbox{span}\,\{a|x\rgl\lgl z|: a\in D_S, \; x, z\in E\}.
\eean
Take $(a,\ol{S}a)\in G(\ol{S})$. There exists $a_n\in D_S$ such that
$a=\lim_n a_n$, $\ol{S}a=\lim_n Sa_n$. Since $E$ is separable, so is $\clk(E)$.
Hence by proposition~1.7.2 in \cite{dix}, there is an approximate identity
$\{p_n\}\in \mbox{span}\,\{|x\rgl\lgl z|:x,z\in E\}$ such that
$\|bp_n-b\|$ converges to zero for all $b\in\clk(E)$, and $\|p_n\|\leq 1$.
Hence $\|a_np_n-a\|\leq \|a_n-a\|\|p_n\|+\|ap_n-a\|$, which implies
$a=\lim_n a_np_n$. Therefore
\bean
\|S(a_np_n)-\ol{S}a\| & \leq &
\|(Sa_n)p_n-(\ol{S}a)p_n\|+\|(\ol{S}a)p_n-\ol{S}a\|\\
 &\leq & \|Sa_n-\ol{S}a\|+\|(\ol{S}a)p_n-\ol{S}a\|,
\eean
and consequently, $\ol{S}a=\lim_n S(a_np_n)$. Since $a_np_n\in
D(\phi_1(\phi_2(S)))$,
$D(\phi_1(\phi_2(S)))$ is a core for $\ol{S}$.

Next, take $(x,\ol{T}x)\in G(\ol{T})$. There exist $\{x_n\}\in D_T$ such that
$x=\lim_n x_n$, and $\ol{T}=\lim_n Tx_n$. Now
$\mbox{span}\,\{\lgl y,z\rgl : y,z \in E\}$
is a dense two-sided ideal in $\ol{\lgl E,E\rgl}$. Again, by proposition~1.7.2
in \cite{dix}, it admits an approximate identity $\{\xi_n\}$ of
$\ol{\lgl E, E\rgl}$ with $\|\xi_n\|\leq 1$. Check that
$x=\lim_n x_n\xi_n$. Since
\bean
D(\phi_2(\phi_1(T))) &=&\mbox{span}\,\{az:a\in D(\phi_1(T)), z\in E\}\\
&=& \mbox{span}\,\{|x\rgl\lgl y|z:x\in D_T, y,z\in E\}\\
&=&\mbox{span}\,\{x\lgl y,z\rgl:x\in D_T, y,z\in E\},
\eean
$x_n\xi_n$ is in $D(\phi_2(\phi_1(T)))$ for all $n$.
Therefore $\ol{T}x=\lim_n Tx_n=\lim_n(Tx_n)\xi_n=\lim_n T(x_n\xi_n)$.
Thus $D(\phi_2(\phi_1(T)))$ is a core for $\ol{T}$.
\qed

Let us call two semiregular operators $T_1$ and $T_2$ equivalent if
their closures are equal. In such a case we will write $T_1\sim T_2$.
Clearly this is an equivalence relation.

\blmma \label{t3.5}
Let $S, S_1, S_2\in SR(\clk(E))$ and $T, T_1, T_2\in SR(E)$. Then one has
\[
\brray{rl}
1.&\phi_1\phi_2(S)\sim S, \quad \phi_2\phi_1(T)\sim T,\\
2.& T_1\sim T_2 \iff \phi_1(T_1)\sim \phi_1(T_2),\\
3.& S_1\sim S_2 \iff \phi_2(S_1)\sim \phi_2(S_2).
\erray
\]
\elmma
\prf Part~1 is a consequence of lemma~\ref{t3.3} and lemma~\ref{t3.4}.

For part~2, assume $T_1\sim T_2$, i.e.\ $\ol{T_1}=\ol{T_2}$. Then
$\phi_1(\ol{T_1})=\phi_1(\ol{T_2})$. By lemma~\ref{t3.2},
$\phi_1(T_1)\seq\phi_1(\ol{T_1})\seq \ol{\phi_1(T_1)}$.
Therefore $\ol{\phi_1(T_1)}=\ol{\phi_1(\ol{T_1})}$.
Similarly $\ol{\phi_1(T_2)}=\ol{\phi_1(\ol{T_2})}$.
Hence $\ol{\phi_1(T_1)}=\ol{\phi_1(T_2)}$. Conversely, if
$\ol{\phi_1(T_1)}=\ol{\phi_1(T_2)}$, one has
$\phi_2(\ol{\phi_1(T_1)})=\phi_2(\ol{\phi_1(T_2)})$. By lemma~\ref{t3.2}
\[
\phi_2\phi_1(T_i)\seq \phi_2(\ol{\phi_1(T_i)})\seq\ol{\phi_2\phi_1(T_i)},
\quad i=1,2.
\]
Hence by part~1, $\ol{T_i}=\ol{\phi_2(\ol{\phi_1(T_i)})}$, $i=1,2$,
which now implies $\ol{T_1}=\ol{T_2}$.

Proof of part~3 is exactly similar.\qed

If we denote by $sr(E)$ (respectively  $sr(\clk(E))$) the space of all
semiregular operators on $E$ (resp.\ $\clk(E)$) modulo the above
equivalence relation, then the above lemma tells us that the maps
$\phi_1:sr(E)\raro sr(\clk(E))$
and
$\phi_2:sr(\clk(E))\raro sr(E)$
are one-one, onto and are inverses of each other. Therefore we can identify
a semiregular operator $T$ on $E$ with its image $\phi_1(T)$ on $\clk(E)$.
Also, from the definitions of $\phi_1$ and $\phi_2$ it follows easily that
if $T$ (respectively $S$) is regular on $E$ (resp.\ $\clk(E)$) with
$z$-transform $z_T$ (resp.\ $z_S$), then $\ol{\phi_1(T)}$
(resp.\ $\ol{\phi_2(S)}$) is regular on $\clk(E)$ (resp.\ $E$) with
the same $z$-transform. Thus the identification of semiregular operators that
we are making is  compatible with the identification
of regular operators on the two spaces that we have already made in the
beginning  of this section using their $z$-transforms.

\newsection{Abelian  $C^*$-algebras} \label{s4}
In this section, we will prove that on  $C^*$-algebras of the form
$C_0(X)$, where $X$ is a locally compact Hausdorff space,
any semiregular operator is given by multiplication by a
continuous function, thereby implying that it is regular.

\bppsn \label{t4.1}
Let $X$ be a locally compact Hausdorff space, and let $T$ be a closed
semiregular operator on $C_0(X)$. Then $T$ is regular.
\eppsn

\prf The proof is quite elementary. The key observation in the proof
is the fact that the Pedersen ideal of $C_0(X)$ is $C_c(X)$, that is,
any dense ideal in $C_0(X)$ contains $C_c(X)$ (see 5.6.3, p-176,\cite{pe}).
 Let $\LMD$ be the set
$\{K\seq X: K\mbox{ compact}\}$ ordered by inclusion.
For each $K\in\LMD$, choose a $f_K\in C_c(X)$ such that
$0\leq f_K(x)\leq 1$ for all $x$, and $f_K(x)=1$ for all $x$ in $K$.
We have already observed that $C_c(X)\seq D_T$, so that each $f_K$ is in the
domain of $T$. Our claim now is that the net $\{Tf_K\}_{K\in\LMD}$
converges pointwise to a continuous function $f$ on $X$. First, let us show
that for any $x\in X$, the net $\{Tf_K(x)\}_{K\in\LMD}$ converges.
Let $K_0=\mbox{supp}\,f_{\{x\}}$. Then for any $K\supseteq K_0$,
$f_Kf_{\{x\}}=f_{\{x\}}$. Therefore $Tf_{\{x\}}=(Tf_K)f_{\{x\}}$.
Evaluating at the point $x$, we get $(Tf_K)(x)=(Tf_{\{x\}})(x)$
whenever $K_0\seq K$. So $\{Tf_K(x)\}_{K\in\LMD}$ converges.
Define $f(x):=\lim_K Tf_K(x)$.  Take any compact subset $K$ of $X$.
Let $S_K$ be the support of $f_K$. Then for any $K_0\supseteq S_K$,
$Tf_{K_0}(x)=Tf_K(x)$ for all $x\in K$. Hence $f(x)=Tf_K(x)$ for all $x\in K$.
Since $Tf_K\in C_0(X)$, $f$ is continuous on $K$. This being true for
any compact subset $K$ of $X$, $f$ is continuous on $X$.

Next observe that $Tg=fg$ for all $g\in C_c(X)$. Indeed, if
$K=\mbox{supp}\,g$, then $g=f_Kg$. Therefore
$Tg=(Tf_K)g$. Since $Tf_K=f$ on $K$, and $g=0$ outside $K$, we have $Tg=fg$.
If we denote by $T_f$ the operator $g\mapsto fg$ on $C_0(X)$ (with
maximal domain), then $T_f$ is regular, $T_f|_{C_c(X)}=T|_{C_c(X)}$.
Let $D=(1+|f|^2)^{1/2}C_c(X)$. It is easy to see that $D$ is dense in $C_0(X)$.
Therefore $C_c(X)=(1+|f|^2)^{-1/2}D$ is a core for $T_f$. So
$T_f=\ol{T_f|_{C_c(X)}}=\ol{T|_{C_c(X)}}$. Since $T$ is closed, this implies
\be \label{e4.1}
T_f\seq T,
\ee
and hence $T^*\seq {T_f}^*$. Since $D_{T^*}$ is a dense ideal
in $C_0(X)$, we get $C_c(X)\seq D_{T^*}$.
So $T^*|_{C_c(X)}={T_f}^*|_{C_c(X)}$. Now, ${T_f}^*$ is just multiplication
by $\ol{f}$ and $C_c(X)$ is a core for ${T_f}^*$. Therefore
${T_f}^*=\ol{{T_f}^*|_{C_c(X)}}=\ol{T^*|_{C_c(X)}}\seq T^*$. Since
$T_f$ is regular, this implies
$T^{**}\seq{T_f}^{**}=T_f$, and consequently, $T\seq T_f$. This, along with
(\ref{e4.1}), implies $T=T_f$. Thus $T$ is regular.\qed

\brmrk \label{t4.3}
{\rm
Let $\cla$ be a unital $C^*$-algebra.
Observe that $C_0(X)\ots\cla$ can be identified with the space
$C_0(X,\cla)$ of $\cla$-valued continuous functions on X that vanish
at infinity, with its usual norm. Notice also that Tietze's extension theorem
continues to hold for $\cla$-valued continuous functions. Using this,
it is not too difficult to show that
if D is a dense right ideal
in  $C_0(X,\cla)$, then it must contain  $C_c(X,\cla)$,
the space of all compactly supported $\cla$-valued functions on $X$.
Having proved this, notice now that the proof of proposition~\ref{t4.1}
remains valid if one replaces $C_0(X)$ by $C_0(X,\cla)$ and
$C_c(X)$ by $C_c(X,\cla)$. Thus proposition~\ref{t4.1} continues to hold
for semiregular operators on $C^*$-algebras of the form $C_0(X)\ots\cla$
as well.
}
\ermrk
\newsection{Subalgebras of $\clb_0(\clh)$} \label{s5}
We will deal with non abelian $C^*$-algebras in this section. The
simplest case of course is the algebra $\clb_0(\clh)$ of compact operators
on a Hilbert space. We have seen in section~\ref{s3} that semiregulars on
$\clb_0(\clh)$ are, up to taking closures, same as semiregular operators
on $\clh$. Therefore it is natural to expect that in this case, all
semiregular operators are regular. The following proposition says
that this is indeed the case.

\bppsn \label{t5.1}
Let $\clh$ be a complex separable Hilbert space. Then any closed semiregular
operator on $\clb_0(\clh)$ is regular.
\eppsn
\prf This is a straightforward consequence of the results in section~3.
Let $T$ be a closed semiregular operator on $\clb_0(\clh)$. Then
$\phi_2(T)$ is a semiregular operator on $\clh$, which simply means that
it is closable and densely defined. Therefore its closure $\ol{\phi_2(T)}$ is
a regular operator on $\clh$. From the remarks following lemma~\ref{t3.5},
$\ol{\phi_1(\ol{\phi_2(T)})}$ is a regular operator on $\clb_0(\clh)$.
By lemma~\ref{t3.3} and~\ref{t3.5}, we now obtain
$\phi_1\phi_2(T)\seq\phi_1\ol{\phi_2(T)}\seq\ol{\phi_1\phi_2(T)}=\ol{T}=T$.
Hence $T=\ol{\phi_1(\ol{\phi_2(T)})}$. But since $\ol{\phi_2(T)}$ is
regular, so is $\ol{\phi_1(\ol{\phi_2(T)})}$. Thus $T$ is regular.\qed

Let us now consider the next simplest class of $C^*$-algebras, namely,
the sub $C^*$-algebras of $\clb_0(\clh)$.

Let $\cla$ be a $C^*$-algebra. Denote by $\hta$ its spectrum, i.e.\ the space
of all
irreducible representations of $\cla$, equipped with its usual
hull-kernel topology. Let $T$ be a semiregular operator
on $\cla$, with domain $D(T)$. Since $D(T)$ is a dense right ideal in $\cla$,
for any $\pi\in\hta$,
$D(T)_\pi:=\{\pi(a):a\in D(T)\}$ is a dense right ideal in the $C^*$-algebra
$\pi(\cla)$. Define an operator $T_\pi$ on $D(T)_\pi$ by the prescription
\be \label{e5.0}
T_\pi \pi(a):=\pi(Ta),\quad a\in D(T).
\ee
To see that this is well defined, notice that if $\pi(a)=\pi(b)$,
then for any $c\in D(T^*)$,
\[
\lgl \pi(Ta)-\pi(Tb),\pi(c)\rgl=\lgl \pi(a)-\pi(b), \pi(T^*c)\rgl=0.
\]
Since $D(T^*)$ is dense in $\cla$, $\pi(Ta)=\pi(Tb)$.  The equality
\[
\lgl \pi(Ta),\pi(c)\rgl=\lgl\pi(a),\pi(T^*c)\rgl\quad
                         \forall a\in D(T),c\in D(T^*)
\]
shows that $T_\pi$ is closable and $(T^*)_\pi\seq (T_\pi)^*$,
therby implying that ${T_\pi}^*$ is densely defined. Thus $T_\pi$ is
a semiregular operator on $\pi(\cla)$.

\blmma \label{t5.2}
If $T$ is regular, then each $\ol{T_\pi}$ is regular.
\elmma
\prf All we need to show is that $\mbox{ran}\,(I+{T_\pi}^*T_\pi)=\pi(\cla)$.
Take any $b\in\cla$. By regularity of $T$, there is an $a\in D(T^*T)\seq D(T)$
such that $(I+T^*T)a=b$. But then
\bean
(I+{T_\pi}^*T_\pi)\pi(a) &=& \pi(a)+{T_\pi}^*T_\pi \pi(a)\\
  &=& \pi(a)+{T_\pi}^*\pi(Ta)\\
  &=&\pi(a)+(T^*)_\pi\pi(Ta)\\
  &=& \pi(a)+ \pi(T^*Ta)\\
  &=& \pi((I+T^*T)a)\\
  &=& \pi(b).
\eean
Thus $\mbox{ran}\,(I+{T_\pi}^*T_\pi)=\pi(\cla)$.\qed

\blmma \label{t5.3}
Let $S$ and $T$ be semiregular operators on $\cla$. Then for each
$\pi\in\hta$,\\[-3ex]
{\rm
\begin{enumerate}\renewcommand{\theenumi}{\arabic{enumi}}
\item $S\seq T\RARO S_\pi\seq T_\pi$,
\item $D(S^*S)_\pi\seq D({S^*}_\pi S_\pi)\seq D({S_\pi}^*S_\pi)$,
\item {\it if $D(S^*S)$ is a core for $S$, then $D(S^*S)_\pi$ is a core for
$S_\pi$.}
\end{enumerate}
}
\elmma
\prf
Proof of the first two parts are trivial. For part~3, take
$(\pi(a), \ol{S_\pi}\pi(a))\in G(\ol{S_\pi})$. Choose $a_0\in D(S)$ such that
\[
\|(\pi(a_0),S_\pi \pi(a_0))- (\pi(a), \ol{S_\pi}\pi(a))\|<\eps.
\]
Since $D(S^*S)$ is a core for $S$, there is an $a_1\in D(S^*S)$ such that
$\|(a_1,Sa_1)-(a_0,Sa_0)\|<\eps$. But then
$\|(\pi(a_1), \pi(Sa_1))- (\pi(a_0), \pi(Sa_0))\|<\eps$,
so that
\[
\|(\pi(a_1),S_\pi \pi(a_1))- (\pi(a), \ol{S_\pi}\pi(a))\|<2\eps.
\]
Thus $D(S^*S)_\pi$ is a core for $S_\pi$.\qed

\bppsn \label{t5.4}
Let $S$ be a closed semiregular operator on $\cla$. If each $\ol{S_\pi}$
is regular, and $D(S^*S)_\pi$ is a core for ${S_\pi}^*S_\pi$ for
all $\pi\in\hta$, then $S$ is regular.
\eppsn
\prf
The given condition implies that any element of the form
$(I+{S_\pi}^*S_\pi)\pi(a)$ can be approximated by an element of
the form $(I+{S_\pi}^*S_\pi)\pi(b)$ where $b\in D(S^*S)$, which, in turn,
implies that $\pi\{(I+S^*S)b:b\in D(S^*S)\}$
is dense in $\pi(\cla)$ for all $\pi\in\hta$. By theorem~\ref{t1.1},
$\mbox{ran}\,(I+S^*S)$ is dense in $\cla$.\qed

\blmma \label{t5.5}
Assume that $\cla$ is separable or GCR, and $\hta$ has discrete topology.
Then for any $a\in\cla$ and $\pi\in\hta$, there is a unique element
$a_{(\pi)}\in\cla$ such that
\[
\pi'(a_{(\pi)})=\cases{0 & if $\pi'\neq\pi$,\cr
                      \pi(a) & if $\pi'=\pi$.}
\]
\elmma
\prf
Let $J=\cap_{\pi'\neq\pi}\mbox{ker\,}\pi'$.
If $J\seq\kernel\pi$, then $\pi$ will belong to the closure of
$\{\pi\}^c$. But $\{\pi\}$ is open. Hence $J\not\seq\kernel\pi$.
Therefore $\pi|_J$ is nonzero, and since $J$ is an ideal,
$\pi|_J$ is actually irreducible. From the assumptions, it follows that
$\cla$ is liminal; therefore we get
$\pi(J)=\clb_0(\clh_\pi)=\pi(\cla)$.
So there is an element $a_{(\pi)}\in J$ such that
$\pi(a_{(\pi)})=\pi(a)$. Obviously $\pi'(a_{(\pi)})=0$ for
all other $\pi'$. Uniqueness is now obvious.%
\qed

\blmma \label{t5.6}
Under the same assumptions as in the previous lemma, we have\\[-3ex]
{\rm
\begin{enumerate}\renewcommand{\theenumi}{\arabic{enumi}}
\item $(ab)_{(\pi)}=a\cdot b_{(\pi)}$,
\item $\pi(a)=\lim_n\pi(a^{(n)}) \RARO a_{(\pi)}=\lim_n a^{(n)}_{(\pi)}$,
\item {\it If $S$ is a closed semiregular operator on $\cla$, then
 $a\in D(S)$ implies $a_{(\pi)}\in D(S)$ and
 $Sa_{(\pi)}=(Sa)_{(\pi)}$ for all $\pi\in\hta$.}
\end{enumerate}
}
\elmma
\prf
Parts 1 and 2 are obvious. For part~3, choose an approximate identity
$\{e^{(n)}\}$ in $\cla$. It is easy to verify that $ae^{(n)}_{(\pi)}$
converges to $a_{(\pi)}$ and $S(ae^{(n)}_{(\pi)})$ converges. The result
follows by closedness of $S$.\qed

\blmma \label{t5.7}
Let $\cla$ be as above, and let $S$ be a closed semiregular operator
on $\cla$. Then\\[-3ex]
{\rm
\begin{enumerate}\renewcommand{\theenumi}{\arabic{enumi}}
\item {\it $S_\pi$ is closed,}
\item ${S_\pi}^*=(S^*)_\pi$, {\it and}
\item $(S^*S)_\pi={S_\pi}^*S_\pi$.
\end{enumerate}
}
\elmma
\prf
1. Take $\pi(a)\in D(\ol{S_\pi})$. There exist $a^{(n)}\in D(S)$
such that $\pi(a^{(n)})$ converges to $\pi(a)$ and $S_\pi\pi(a^{(n)})$
converges. By the previous lemma, $a^{(n)}_{(\pi)}$ converges
to $a_{(\pi)}$, and $(Sa^{(n)})_{(\pi)}$ converges. By closedness
of $S$, we conclude that $a_{(\pi)}\in D(S)$ and
$Sa_{(\pi)}=\lim (Sa^{(n)})_{(\pi)}$. Therefore
$\pi(a)=\pi(a_{(\pi)})\in D(S)_\pi=D(S_\pi)$.

2. We only have to show the inclusion ${S_\pi}^*\seq (S^*)_\pi$.
Take $\pi(a)\in D({S_\pi}^*)$. For any $b\in D(S)$ and
$\pi'\in\hta$, $\pi'\neq\pi$, we have
$\pi'(\lgl Sb,a_{(\pi)}\rgl)=\pi'(Sb)^*\pi'(a_{(\pi)})=0$,
and
\bean
\pi(\lgl Sb,a_{(\pi)}\rgl) &=& \pi(Sb)^*\pi(a_{(\pi)})\\
 &=&\lgl S_\pi \pi(b), \pi(a)\rgl\\
 &=&\lgl \pi(b), {S_\pi}^*\pi(a)\rgl.
\eean
Since ${S_\pi}^*\pi(a)\in \pi(\cla)$, ${S_\pi}^*\pi(a)=\pi(c)$ for some
$c\in\cla$. Hence
$\pi(\lgl Sb, a_{(\pi)}\rgl)=\lgl\pi(b),\pi(c)\rgl=\pi(\lgl b,c\rgl)$.
Therefore $\pi'(\lgl Sb, a_{(\pi)}\rgl)=\pi'(\lgl b,c_{(\pi)}\rgl)$
for all $\pi'$, i.e.\ $\lgl Sb,a_{(\pi)}\rgl =\lgl b,c_{(\pi)}\rgl$.
This implies $a_{(\pi)}\in D(S^*)$ and consequently
$\pi(a)=\pi(a_{(\pi)})\in D(S^*)_\pi$.

3. Again, the inclusion $D(S^*S)_\pi\seq D({S_\pi}^*S_\pi)$
is obvious. To show the reverse inclusion, take
$\pi(a)\in D({S_\pi}^*S_\pi)$. Then $\pi(a)\in D(S_\pi)=D(S)_\pi$,
and $S_\pi \pi(a)\in D({S_\pi}^*)=D(S^*)_\pi$.
This means $a_{(\pi)}\in D(S)$, and $(Sa)_{(\pi)}=Sa_{(\pi)}\in D(S^*)$,
i.e.\ $a_{(\pi)}\in D(S^*S)$. Hence $\pi(a)=\pi(a_{(\pi)})\in D(S^*S)_\pi$.%
\qed

As a consequence of the above results, we now obtain the following.

\bthm \label{t5.8}
Let $\clh$ be a complex separable Hilbert space. Any closed semiregular
operator on a $C^*$-subalgebra of $\clb_0(\clh)$ is regular.
\ethm
\prf
Follows from proposition~\ref{t5.4} and lemma~\ref{t5.7}.\qed

Let $G$ be a compact quantum group (\cite{wo2}), and let $\whtG$  denote
its Pontryagin dual. Then the $C^*$-algebra $C_0(\whtG)$ has discrete spectrum.
Therefore in this context, we can rephrase the previous result as follows.

\bppsn \label{t5.9}
Let $G$ be a compact quantum group. Any closed semiregular operator
on $C_0(\whtG)$ is regular.
\eppsn

\newsection{Extensions  of semiregular operators} \label{s6}
We will be concerned with more general classes of $C^*$-algebras
in this section. As the example in section~\ref{s2} suggests,
we can not possibly expect results like theorem~\ref{t5.8} to
hold once we go beyond subalgebras of $\clb_0(\clh)$. However
in many cases, it is possible to get regular extensions of semiregular
operators. But before we go to extensions of semiregular operators,
let us find out how a regular operator is related to its extensions
and restrictions.

\bppsn \label{t6.1}
Let $T$ be a regular operator on a $C^*$-algebra $\cla$ with $z$-transform $z$.
Let $u$ be an isometry in $M(\cla)$ obeying the following condition:
\be \label{e6.1}
(u^*(I-z^*z)u)^{1/2}=(I-z^*z)^{1/2}u.
\ee
Then $z_S:=zu$ is the $z$-transform of a regular restriction $S$ of $T$.
\eppsn
\prf
Since $z$ is the $z$-transform of a regular operator,
$(I-z^*z)^{1/2}\cla$ is dense in $\cla$, and $u$ is an isometry.
So $u^*(I-z^*z)^{1/2}\cla$ is also dense in $\cla$. Now
from the given conditions, we get
\[
(I-z^*z)^{1/2}u=u^*(I-z^*z)^{1/2}.
\]
Therefore $(I-z^*z)^{1/2}u\cla$ is dense in $\cla$. It follows then that
$(u^*(I-z^*z)u)\cla=(I-{z_S}^*z_S)\cla$ is dense in $\cla$. This means
$(I-{z_S}^*z_S)^{1/2}\cla$ is also dense in $\cla$. Clearly $\|z_S\|\leq 1$.
Hence there is a unique regular operator $S$ whose $z$-transform is $z_S$.

To show that $S$ is a restriction of $T$, it is enough to prove that $S=T$
on $(I-{z_S}^*z_S)\cla$, since this is a core for $S$. We will prove that
(i)~$(I-{z_S}^*z_S)\cla\seq (I-z^*z)^{1/2}\cla$, and
(ii)~if $(I-{z_S}^*z_S)a=(I-z^*z)^{1/2}b$, then
$z_S(I-{z_S}^*z_S)^{1/2}a=zb$.
(i) is a direct consequence of (\ref{e6.1}). For (ii), assume that
$(I-{z_S}^*z_S)a=(I-z^*z)^{1/2}b$. This means
$(u^*(I-z^*z)u)a=(I-z^*z)^{1/2}b$, which, together with (\ref{e6.1}) and
injectivity of the operator $(I-z^*z)^{-1/2}$,
implies that $u(I-z^*z)^{1/2}ua=b$. Therefore
$zb=zu(I-z^*z)^{1/2}ua=z_S(I-z^*z)^{1/2}ua=z_S(I-{z_S}^*z_S)^{1/2}a$.
\qed
\bppsn \label{t6.2}
Let $T$ be a regular operator on $\cla$ with $z$-transform $z$. A regular
operator $S$ with $z$-transform $z_S$ is a restriction of $T$ if and only if
$z_S =zu$ for some isometry $u$ in $M(\cla)$ obeying equation~(\ref{e6.1}).
\eppsn
\prf
We have seen that if $u$ is such an isometry, then $z_S:=zu$ defines
a regular restriction of $T$. Now conversely, suppose $S$ is a regular
restriction of $T$. Since $D(S)\seq D(T)$, we have
$(I-{z_S}^*z_S)^{1/2}\cla\seq (I-z^*z)^{1/2}\cla$. Therefore the
operator $w_{S,T}:= (I-z^*z)^{-1/2}(I-{z_S}^*z_S)^{1/2}$ is everywhere
defined on $\cla$. $(I-{z_S}^*z_S)^{1/2}$ is bounded,
$(I-z^*z)^{-1/2}$ is closed; so $w_{S,T}$ is also closed. Hence
$w_{S,T}$ is bounded. Observe next that for any $a,b\in\cla$,
\[
\lgl w_{S,T}a,(I-z^*z)^{1/2}b\rgl=\lgl (I-{z_S}^*z_S)^{1/2}a, b\rgl
  =\lgl a, (I-{z_S}^*z_S)^{1/2}b\rgl.
\]
This implies $D(T)\seq D({w_{S,T}}^*)$, i.e.\ ${w_{S,T}}^*$ is
densely defined. Together with the boundedness of $w_{S,T}$, this
means it is adjointable.

For any $a\in\cla$,
$(I-{z_S}^*z_S)^{1/2}a=(I-z^*z)^{1/2}w_{S,T}a\in D(S)\seq D(T)$.
Since $S=T$ on $D(S)$, we have $z_Sa=zw_{S,T}a$ for all $a\in\cla$.
So $z_S=zw_{S,T}$. Next, take any $a,b\in\cla$. Then
\bean
\lgl {w_{S,T}}^*(I-z^*z)w_{S,T}a,b\rgl &=&
                   \lgl (I-z^*z)^{1/2}w_{S,T}a, (I-z^*z)^{1/2}w_{S,T}b\rgl\\
 & =&\lgl (I-{z_S}^*z_S)a,b\rgl.
\eean
Therefore ${w_{S,T}}^*(I-z^*z)w_{S,T}=I-{z_S}^*z_S$.
It now follows immediately that ${w_{S,T}}^*w_{S,T}=I$, and
$\bigl({w_{S,T}}^*(I-z^*z)w_{S,T}\bigr)^{1/2}=(I-z^*z)^{1/2}w_{S,T}$.
\qed

\bppsn \label{t6.3}
Let $T$ be a regular operator on $\cla$ with $z$-transform $z$. A regular
operator $S$ with $z$-transform $z_S$ is an extension of $T$ if and only if
$z_S=uz$ for some coisometry $u$ satisfying the following equation:
\be\label{e6.2}
(u(I-zz^*)u^*)^{1/2}=u(I-zz^*)^{1/2}.
\ee
\eppsn
\prf
Use the previous proposition and the fact that $S$ is an extension
of $T$ if and only if $S^*$ is a restriction of $T^*$.
\qed

Let $S$ be a closed semiregular operator on a $C^*$-algebra $\cla$. For any
$\pi\in\hta$, if we define $S_\pi$ by equation~(\ref{e5.0}), with $S$
replacing $T$, then $\ol{S_\pi}$ is a closed semiregular operator
on $\pi(\cla)$. Construct an operator $\tlds$ on $\cla$ as follows:
\be \label{e6.3}
\brray{c}
D(\tlds)=\{a\in\cla:\pi(a)\in D(\ol{S_\pi})\:\forall\, \pi\in\hta,
 \; \exists\,b\in\cla\,\ni\,\pi(b)=\ol{S_\pi}\pi(a)\:\forall \pi\in\hta\},\\
\tlds a=b.
\erray
\ee

\blmma \label{t6.4}
$\tlds$ is a closed semiregular extension of $S$.
\elmma
\prf
{}From the definition of $D(\tlds)$ and the fact that $\hta$ separates
points of $\cla$, it follows that $\tlds$ is well defined.
To show it is closed, take $a_n\in D(\tlds)$ such that $a_n$
converges to $c$ and $\tlds a_n$ converges
to $d$. Then for each $\pi\in\hta$, $\pi(a_n)$ converges to $\pi(a)$ and
$\pi(\tlds a_n)$ converges to $\pi(d)$. But
$\pi(\tlds a_n)=\ol{S_\pi}\pi(a_n)$.
Therefore by closedness of $\ol{S_\pi}$, $\pi(c)\in D(\ol{S_\pi})$
and $\ol{S_\pi}\pi(c)=\pi(d)$. This means $c$ is in $D(\tlds)$ and
$\tlds c =d$, i.e.\ $\tlds$ is closed. The inclusion $S\seq \tlds$ is
obvious. So $\tlds$ is densely defined.

Using the same arguments for $S^*$, we see that $\widetilde{(S^*)}$ is
densely defined. It is routine to check that for any $a\in D(\tlds)$
and $b\in D(\widetilde{(S^*)})$, one has
$\lgl \tlds a, b\rgl = \lgl a, \widetilde{(S^*)} b\rgl$,
so that $\widetilde{(S^*)}\seq (\tilde{S})^*$.
Hence $(\tilde{S})^*$ is densely defined.

Thus $\tlds$ is a closed semiregular extension of $S$.
\qed

Next we compute the adjoint of $\tlds$.

\blmma\label{t6.5}
Let $\tlds$ be as above. Then $(\tilde{S})^*=\widetilde{(S^*)}$.
\elmma
\prf
Define an operator $S_*$ by the prescription
\[
D(S_*)=\{a\in\cla: \pi(a)\in D({S_\pi}^*), \; \exists\, b\in\cla \ni
                      \pi(b)={S_\pi}^* \pi(a) \:\forall \pi\in\hta \},
\]
\[
S_*a=b.
\]
We claim that
\[
 \widetilde{(S^*)}\seq S_*\seq (\tlds)^*.
\]
The first inclusion follows from the inclusion
$(S^*)_\pi\seq {S_\pi}^*$. The second inclusion follows from the
observation that for any $a\in D(\tlds)$ and $b\in D(S_*)$, one has
$\pi(\lgl \tlds a,b\rgl)=\pi(\lgl a, S_* b\rgl)$ for all $\pi\in\hta$.
Now, we also have, by the previous lemma,
$S\seq \tlds$ and $S^*\seq \widetilde{(S^*)}$. Thus we have the chain of
inclusions
\[
\widetilde{(S^*)}\seq S_*\seq (\tlds)^*\seq S^*\seq \widetilde{(S^*)},
\]
which proves the lemma.
\qed

\blmma \label{t6.6}
If $S$ is regular, then so is $\tlds$, and $S=\tlds$.
\elmma
\prf
Let $\tau$ denote the operator
$a\oplus b\mapsto b\oplus (-a)$ on $\cla\oplus\cla$.
Then for any operator $T$, $\tau(G(T^*))\seq G(T)^\perp$.
Hence $\tau\bigl(G((\tlds)^*)\bigr)\seq G(\tlds)^\perp$.
Now $G(S)\seq G(\tlds)$, so that $G(\tlds)^\perp\seq G(S)^\perp$.
Since $S$ is regular, we have, using the previous lemma,
\[
\tau(G((\tlds)^*))\seq G(\tlds)^\perp\seq G(S)^\perp
      =\tau(G(S^*))\seq \tau(G(\widetilde{S^*}))=\tau(G((\tlds)^*)).
\]
So we actually have equality everywhere. Hence
$\cla\oplus\cla =G(S)\oplus G(S)^\perp \seq G(\tlds)\oplus G(\tlds)^\perp$.
This implies $G(\tlds)\oplus G(\tlds)^\perp =\cla\oplus\cla$.
Consequently $\tlds$ is regular, and $S=S^{**}=(\tlds)^{**}=\tlds$.
\qed

Now we are ready for the main result in this section. Let $S$ be a
semiregular operator on a liminal $C^*$-algebra $\cla$. Since
$\pi(\cla)=\clb_0(\clh_\pi)$, by proposition~\ref{t5.1},
each $\ol{S_\pi}$ is regular. Let $z_\pi$ be the corresponding
$z$-transform.

\bthm \label{t6.7}
Let $S$ and $z_\pi$ be as above. Suppose $\{u_\pi\}_{\pi\in\hta}$
is a family of coisometries, where
$u_\pi\in \clb(\clh_\pi)$, such that
\[
(u_\pi(I-z_\pi{z_\pi}^*){u_\pi}^*)^{1/2}=u_\pi(I-z_\pi{z_\pi}^*)^{1/2},
\]
and there exists an element $z\in M(\cla)$ for which
$\pi(z)=u_\pi z_\pi$ for all $\pi\in\hta$. Then $S$ admits a
regular extension.
\ethm
\prf
Let us first of all show that $z$ is indeed the $z$-transform of a
regular operator. Clearly $\|z\|\leq 1$. By proposition~\ref{t6.3},
each $\pi(z)$ is the $z$-transform of a regular operator on $\pi(\cla)$.
Therefore $(I-\pi(z)^*\pi(z))^{1/2}\clh_\pi$ is dense in $\clh_\pi$.
Since $\pi((I-z^*z)^{1/2})=(I-\pi(z)^*\pi(z))^{1/2}$, it follows from
proposition~2.5 in \cite{wo1} that $(I-z^*z)^{1/2}\cla$ is dense
in $\cla$. Thus $z$ is the $z$-transform of a regular operator, say, $T$.

Our next job is to show that $T$ is an extension of $S$. It is easy to see
that $T_\pi$, defined by (\ref{e5.0}) on $\pi(\cla)$, is a regular operator
with $z$-transform $\pi(z)$. Proposition~\ref{t6.3} tells us that
$\ol{S_\pi}\seq T_\pi$. From the definition of $\tlds$ it follows that
$\tlds\seq \widetilde{T}$. But by lemma~\ref{t6.4}, $S\seq \tlds$ and by
lemma~\ref{t6.6}, $T=\widetilde{T}$. Therefore $T$ is an extension of $S$.
\qed

\brmrk \label{t6.8}
{\rm
{}From the proof of the above theorem, it is clear that if $S$ is a
semiregular operator on any $C^*$-algebra $\cla$ (not necessarily
liminal) such that the closure of each fibre
$S_\pi$ is regular with $z$-transform $z_\pi$, and there is one single
element $z$ in $M(\cla)$ such that $\pi(z)=z_\pi$ for all $\pi$,
then $\tlds$ is a regular operator.
}
\ermrk

It is now easy to see why the example in section~\ref{s2} fails
to be regular. Each of the fibres $\bbt_\pi$'s is regular, acting on
the same Hilbert space $L_2(0,1)$. But while all the $\bbt_\pi$'s
are equal for $\pi>0$, $\bbt_0$ is different. The same is therefore
true for their $z$-transforms $z_\pi$'s. Hence clearly there can not
be any element in $\cll(E)$ (which are precisely the
$\clb(L_2(0,1))$-valued functions on $[0,1]$ that are both strong and
strong${}^*$-continuous) whose $\pi$-image is the
$z$-transform of $\bbt_\pi$ for all $\pi$.

As a consequence of theorem~\ref{t6.7} and remark~\ref{t6.8}, we
now have the following proposition.

\bppsn \label{t6.9}
Let $S$ be a semiregular operator on a $C^*$-algebra $\cla$. Suppose there
exists a $\pi_0\in\hta$, a regular operator $\bbt$ on $\pi_0(\cla)$,
and a family $\{U_{\pi}\}_{\pi\in\hta}$ of unitary operators
\[
U_{\pi}:\clh_{\pi_0}\mapsto\clh_{\pi},\quad \pi \in \hta,
\]
satisfying the following conditions:\\[-3ex]
{\rm
\begin{enumerate}\renewcommand{\theenumi}{\roman{enumi}}
\item {\it $U_{\pi_0}=I$},
\item {\it $S_\pi\seq U_{\pi}\bbt \,{U_\pi}^*$ for all $\pi$}.
\item {\it for any $a\in \clb(\clh_{\pi_0})$, there is an element $b\in M(\cla)$
    such that $\pi(b)=U_{\pi}a\,{U_\pi}^*$ for all $\pi\in\hta$},
\end{enumerate}
}
Then $S$ admits a regular extension.
\eppsn
\prf
Define operators $T_\pi$ on $\pi(\cla)$ as follows:
\[
D(T_\pi)=U_{\pi}D(\bbt){U_\pi}^*,\quad T_\pi= U_{\pi}\bbt\,{U_\pi}^*.
\]
It is routine to verify that   
each $T_\pi$ is regular,
\[
D({T_\pi}^*)=U_{\pi}D(\bbt^*){U_\pi}^*,\quad
          {T_\pi}^*= U_{\pi}\bbt^*\,{U_\pi}^*,
\]
and if $w$ is the $z$-transform of $\bbt$, then the $z$-transform
of $T_\pi$ is $U_\pi w\,{U_\pi}^*$.
Condition~(iii) now ensures the existence of an element $z$ in $M(\cla)$
such that $\pi(z)=U_\pi w\,{U_\pi}^*$ for all $\pi$. By remark~\ref{t6.8}
above, it follows that $\widetilde{T}$ constructed out of these
$T_\pi$'s by the prescription (\ref{e6.3}) (with $T$ replacing $S$)
is regular. Also, a direct consequence of condition~(ii) above
and the definition of $\widetilde{T}$ is that $S\seq \widetilde{T}$.
\qed

As an immediate corollary of the above proposition,
we can deduce the following.
\bcrlre \label{t6.10}
Let $E=C[0,1]\ots\clh$, $\bbt$ a closed operator on $\clh$, and
$\{U_\pi\}_\pi$ a family of strongly continuous unitaries on $\clh$.
Let $T_\pi={U_\pi}\bbt\,{U_\pi}^*$. Let
\[
D(T)=\{f\in E: f_\pi\in D(T_\pi) \;\forall \pi,
        \pi\mapsto T_\pi f_\pi \mbox{ continuous}\},
\]
\[ (Tf)_\pi=T_\pi f_\pi \]
Then $T$ is a regular operator on $E$.
\ecrlre
\prf
Here $\clk(E)=C[0,1]\ots\clb_0(\clh)$ is the relevant $C^*$-algebra. All
we need to check is that the condition~(iii) in the forgoing
proposition is fulfilled. For any $u\in\clh$, $\pi\mapsto U_\pi u$
is continuous. Hence for any finite-rank operator $S$ on $\clh$,
the function $\pi\mapsto U_\pi S{U_\pi}^*$ is continuous in
the norm topology. By approximating a compact operator by finite-rank
operators, one can show that $\pi\mapsto U_\pi S{U_\pi}^*$ is norm
continuous for compact $S$ also. Next, take any $S\in\clb(\clh)$.
{}From the strong continuity of $\{U_\pi\}$ and $\{{U_\pi}^*\}$,
it follows that for any $\clb_0(\clh)$-valued norm continuous
function $\pi\mapsto R_\pi$, the maps $\pi\mapsto U_\pi S{U_\pi}^*R_\pi$
and $\pi\mapsto R_\pi U_\pi S{U_\pi}^*$ are both norm continuous.
This implies that the function $\pi\mapsto U_\pi S{U_\pi}^*$
is an element of the multiplier algebra $M(C[0,1]\ots\clb_0(\clh))$.
If we call it $b$, then $\pi(b)=U_\pi S{U_\pi}^*$. The third condition
in proposition~\ref{t6.9} is thus satisfied. So $T$ is regular.
\qed

We shall now apply the above result to a specific example.

\bxmpl \label{t6.11}
{\rm
Let $E=C[0,1]\ots L_2(0,1)$, $D$ be as in section~2 
and let $T$ be an operator on $E$ defined
as follows:
\[
D(T)=\{f\in E: f_\pi\in D, f_\pi(1)=e^{i\pi}f_\pi(0),
          \pi\mapsto {f_\pi}'\mbox{ continuous}\},
\]
\[ (Tf)_\pi={f_\pi}'+i\pi f_\pi.\]
Then $T$ is regular.
}
\exmpl
\prf
Just take $U_\pi$ to be multiplication by the function
$x\mapsto e^{i\pi x}$, $\bbt$ to be the operator $T_0$
in section~\ref{s2}, and apply the previous result.
\qed

More generally, Let $g\in E$ obey the following properties:
(i)~$g_\pi\equiv g(\pi,\cdot)$ is absolutely continuous,
(ii)~${g_\pi}'\in L_2(0,1)$,
(iii)~$g(0,x)=0$ for all $x$, and
(iv)~$\pi\mapsto {g_\pi}'$ is continuous.
Then the operator $T$ on $E$ given by the following prescription
is regular:
\[
D(T)=\left\{f\in E: f_\pi\in D,
   f_\pi(1)=\exp\Bigl(i(g(\pi,1)-g(\pi,0))\Bigr)f_\pi(0),
          \pi\mapsto {f_\pi}'\mbox{ continuous}\right\},
\]
\[ (Tf)_\pi={f_\pi}'-i{g_\pi}' f_\pi.\]
In this case one has to take $U_\pi$ to be multiplication by
the function $\exp(ig(\pi,\cdot))$.

{\footnotesize {\bf Acknowledgements.}  %
I would like to thank Debashish Goswami \vspace{-.8ex} for many
useful conversations
and also for initiating the seminar series that eventually led
to this paper.
}


\begin{thebibliography}{99}
\bibitem{bj} Baaj, S. \& Julg, P. : Th\'{e}orie bivariante de
   Kasparov et  op\'{e}rateur non born\'{e}s dans les $C^*$-modules
   hilbertiens, {\em C.\ R.\ Acad.\ Sc.},
     Paris, Series I, 296(1983), 875--878.
\bibitem{dix} Dixmier, J. : {\sl $C^*$-Algebras}, North-Holland, 1977.
\bibitem{hil} Hilsum, M. : Fonctorialit\'{e} en K-th\'{e}ory bivariante
    pour les vari\'{e}t\'{e}s lipschitziennes,
     {\em $K$-Theory}, 3(1989), 401--440.
\bibitem{kap} Kaplansky, I. : Modules over $C^*$-algebras, {\em Amer.\
J.\ Math.}, 75(1953), 839--858.
\bibitem{k} Kasparov, G. G. : Hilbert $C^{*}$-modules --- theorems of
  Stinespring and Voiculescu, {\em J. Operator Theory}, 4(1980), 133--150.
\bibitem{la} Lance, E. C. : {\sl Hilbert $C^*$-modules - A Toolkit for Operator
    Algebraists}, Cambridge University Press, 1995.
\bibitem{pa} Paschke, W. : Inner product modules over $B^*$-algebras, {\em
  Trans.\ Amer.\ Math.\ Soc.}, 182(1973), 443--468.
\bibitem{pe} Pedersen, G. K. : {\sl $C^*$-algebras and their Automorphism
  Groups}, Academic Press, 1979.
\bibitem{rie} Rieffel, P. : Induced representations of $C^*$-algebras,
  {\em Adv.\ in Math.}, 13(1974), 176--257.
\bibitem{wo2} Woronowicz, S. L. : Compact matrix pseudogroups, {\em Comm.\
Math.\   Phys.}, 111(1987), 613--665.
\bibitem{wo1} Woronowicz, S.L. : Unbounded Elements Affiliated With
  $C^*$-algebras and Noncompact Quantum Groups,
    {\em Comm. Math. Phys.}, 136(1991), 399--432.
\end{thebibliography}
\end{document}